\documentclass{slides}
\usepackage{amssymb,latexsym,amsmath,amsthm,graphics,lattice}
\usepackage{graphics}
\usepackage{amscd}

\newtheorem{theorem}{Theorem}
\newtheorem{lemma}{Lemma}

\newcommand{\cpe}{con\-gruence-pre\-serv\-ing extension\xspace}

\newcommand{\Ji}[1]{\tup{J}(#1)}

\newcommand{\R}{\D R}

\newtheorem{problem}{Problem}

\DeclareMathOperator{\Con}{Con}

\newcommand{\ba}{\mbf{a}}
\newcommand{\bb}{\mbf{b}}

\newcommand{\Bo}{B}

\newcommand{\bu}{\mbf{u}}

\newcommand{\prodm}[2]{\prod(\,#1\mid#2\,)}

\begin{document}

\title{Isoform lattices} 
\author{G. Gr\"atzer, E.\,T. Schmidt, and R.\,W. Quackenbush\\
Boulder CO}
\date{}
\maketitle 
 
\begin{slide}
\centerline{\tbf{\large Results}}

Let $L$ be a lattice. We call a congruence relation $\gQ$ of $L$ \emph{isoform}, if any
two congruence classes of $\gQ$ are isomorphic (as lattices). Let us call the lattice
$L$ \emph{isoform}, if all congruences of $L$ are isoform.

G. Gr\"atzer and E.\,T. Schmidt:

\vspace{20pt}

\begin{theorem}\label{T:represent}
Every finite distributive lattice $D$ can be represented as the congruence lattice of a
finite isoform lattice $L$. 
\end{theorem}

G. Gr\"atzer, E.\,T. Schmidt, and R.\,W. Quackenbush:

\vspace{20pt}

\begin{theorem}\label{T:main}
Every finite lattice $K$ has a \cpe to a finite isoform lattice $L$.
\end{theorem}

\end{slide}

\begin{slide}
\centerline{\tbf{\large Background}}

At the beginning of time \dots

\vspace{20pt}

\tbf{Dilworth's theorem.}
Every finite distributive lattice $D$ can be represented as the congruence lattice of a
finite lattice $L$.

\medskip
We want:

Every finite distributive lattice $D$ can be represented as the congruence lattice of a
\tbf{nice} finite lattice $L$.

Such a result is called a \emph{representation theorem}.

Example (G. Gr\"atzer and E. T. Schmidt, 1962):

\vspace{20pt}

\begin{theorem}\label{T:seccomp}
Every finite distributive lattice $D$ can be represented as the congruence lattice of a
sectionally complemented finite lattice $L$.
\end{theorem}

\end{slide} 

\begin{slide}
\centerline{\tbf{\large Topic}}

Why could not lattices be more like groups and rings?

We want representation theorems, where 
congruences behave like in groups and rings.

We shall look at lattices that are

regular

uniform

isoform

\end{slide} 

\begin{slide}
\centerline{\tbf{\large Regular lattices}}

Let $L$ be a lattice. We call a congruence relation $\gQ$ of $L$ \emph{regular}, if any
congruence class of $\gQ$ determines the congruence. Let us call the lattice $L$
\emph{regular}, if all congruences of $L$ are regular. 

Sectionally complemented lattices are regular, so we already have a representation
theorem for finite lattices.
\end{slide} 

\begin{slide}
\centerline{\tbf{\large Uniform lattices}}

Let $L$ be a lattice. We call a congruence relation $\gQ$ of $L$ \emph{uniform}, if any
two congruence classes of $\gQ$ are of the same size. Let us call the
lattice $L$ \emph{uniform}, if all congruences of $L$ are uniform. The following result
was proved by G.~Gr\"atzer, E.\,T. Schmidt, and K. Thomsen, 2003:

\vspace{20pt}

\begin{theorem}\label{T:Uniform}
Every finite distributive lattice $D$ can be represented as the congruence lattice of a
finite uniform lattice $L$. 
\end{theorem}

A uniform lattice is always regular, so the lattices of Theorem~\ref{T:Uniform} are
also regular. 

\end{slide} 

\begin{slide}
Here is the result of the construction for~$D = C_4$.
\vspace{40pt}

\centerline{\scalebox{1.5}{\includegraphics{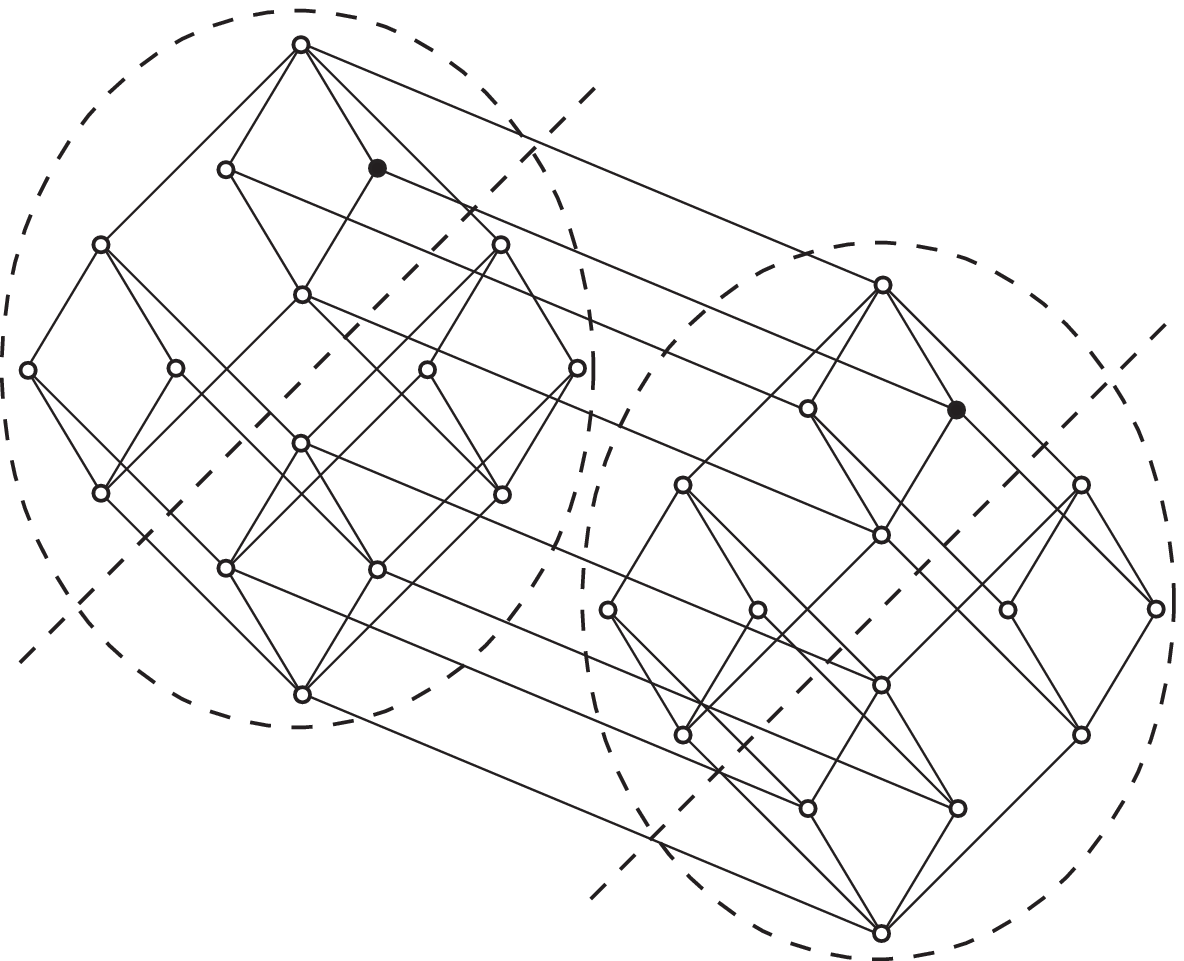}}}

\end{slide}

\begin{slide}
Pruning:

Let $A$ and $B$ be bounded lattices. Redefine $\leq$ on 
$A \times B$ ($A^- = A - \set{0, 1}$):
\[
  \leq_N{}  = {}\leq - \setm{\vv<u, v>}{u,\ v \in A^- \times B,\ u_B \neq v_B}.
\]
$N(A,B) = \vv<A \times B, \leq_N>$. The diagram of $N(B_2,B_1)$:

\vspace{40pt}

\centerline{\includegraphics{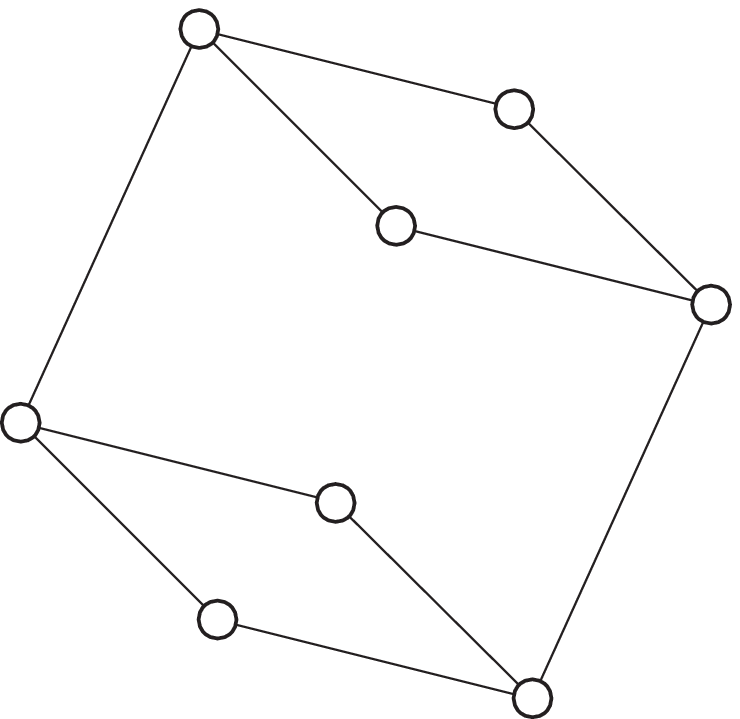}}
\end{slide}

\begin{slide}
The diagram of $N(B_2,B_2)$:

\vspace{40pt}

\centerline{\includegraphics{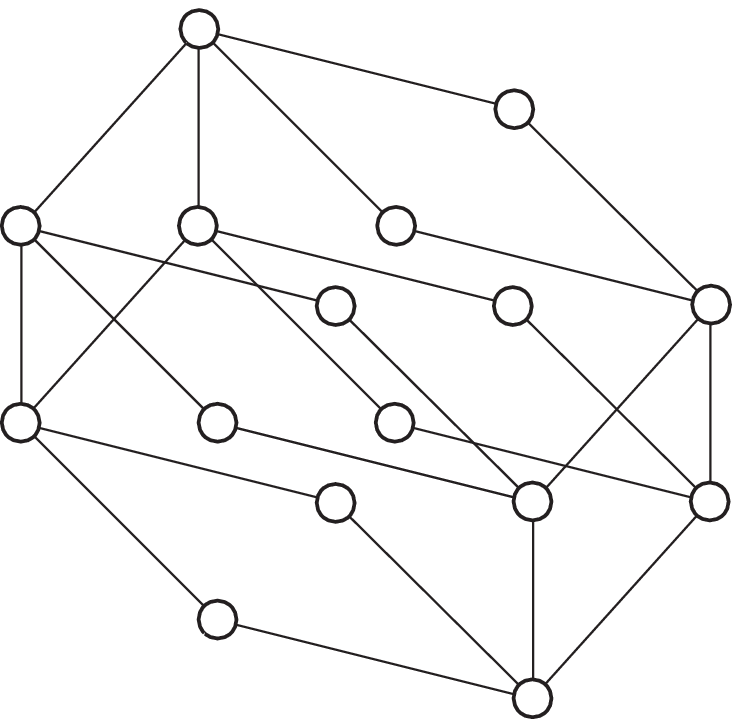}}
\end{slide}

\begin{slide}
\begin{lemma}\label{L:NAB}

$N(A,B)$ is partially ordered by $\leq_N$, in fact, $N(A,B)$ is a lattice. The meet and
join in $N(A,B)$ of $\leq_N$-incomparable elements can be computed using the formulas: 
\begin{align*}
   u \mm_N v &=
     \begin{cases}
        \vv<0,u_B \mm v_B>, &\text{if $u \mm_\times v \in A^- \times B$, $u_B \neq
v_B$};\\
      u \mm_{\times} v, &\text{otherwise.}
     \end{cases}\\
   u \jj_N v &=
      \begin{cases}
          \vv<1,u_B \jj v_B>, &\text{if $u \jj_\times v \in A^- \times B$,  $u_B \neq
v_B$};\\
          u {}\jj_{\times}{} v, &\text{otherwise.}
       \end{cases}
\end{align*}
\end{lemma}

\vspace{40pt}

\begin{lemma}\label{L:cons}
$\Con N(A, B)$ is isomorphic to $\Con A$ with a new zero added.   
\end{lemma}

\end{slide}

\begin{slide}

\begin{theorem}\label{T:uniform}
For any finite distributive lattice $D$, there exists a finite uniform lattice $L$ such
that the congruence lattice of $L$ is isomorphic to $D$, and $L$ satisfies the following
two properties:

(P) Every join-irreducible congruence of $L$ is of the form $\gQ(0,p)$, for a
suitable atom $p$ of $L$. 

(Q) If $\gQ_1$, $\gQ_2$, \dots, $\gQ_n \in \Ji{\Con L}$ are pairwise
incomparable, then $L$  contains atoms $p_1$, $p_2$, \dots, $p_n$ that generate an ideal
isomorphic to $B_n$ and satisfy ${\gQ}_i = \gQ(0,p_i)$, for all $i \leq n$.

\end{theorem}
\end{slide}

\begin{slide}
Let $D$ be a finite distributive lattice with $n$ join-irreducible elements.

If $n = 1$, then $D \iso \Bo_1$.

Let us assume that, for all finite distributive lattices with fewer than $n$
join-irreducible elements, there exists a lattice $L$ satisfying Theorem~\ref{T:uniform}
and properties \tup{(P)} and \tup{(Q)}.

Let $q$ be a minimal element of $\Ji{D}$ and let $q_1$, \dots, $q_k$  ($k \geq 0$) list all
upper covers of $q$ in $\Ji{D}$. Let $D_1$ be a distributive lattice with $\Ji{D_1} = \Ji{D} -
\set{q}$. By the inductive assumption, there exists a lattice $L_1$ satisfying $\Con L_1
\iso D_1$ and \tup{(P)} and \tup{(Q)}.  

If $k = 0$, then $L = \Bo_1 \times L_1$. So we assume that $1 \leq k$.

The congruences of $L_1$ corresponding to the $q_i$'s are pairwise incomparable and
therefore can be written in the form $\gQ(0, p_i)$ and the $p_i$'s generate an ideal $I_1$
isomorphic to $\Bo_k$.
\end{slide}

\begin{slide}
\centerline{\includegraphics{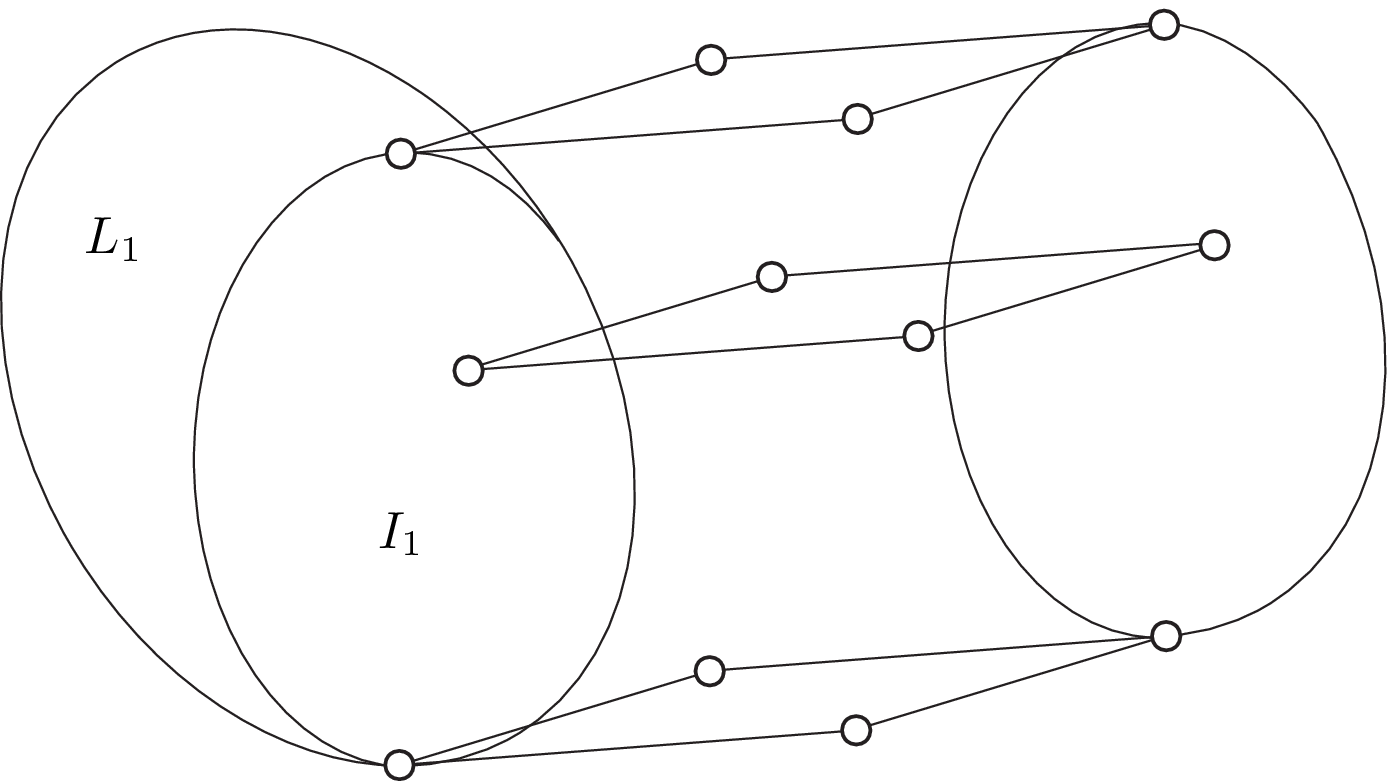}}

\end{slide}

\begin{slide}
\centerline{\tbf{\large Isoform lattices}}

We go beyond groups and rings:

Let $L$ be a lattice. We call a congruence relation $\gQ$ of $L$ \emph{isoform}, if any
two congruence classes of $\gQ$ are isomorphic. Let us call the lattice
$L$ \emph{isoform}, if all congruences of $L$ are isoform.

G. Gr\"atzer and E.\,T. Schmidt:

\vspace{20pt}

\begin{theorem}\label{T:Isoformlect}
Every finite distributive lattice $D$ can be represented as the congruence lattice of a
finite isoform lattice $L$. 
\end{theorem}

Since isomorphic lattices are of the same size, Theorem~\ref{T:Isoformlect} is a
stronger version of Theorem~\ref{T:Uniform}. 

The previous figure shows that a uniform lattice is not necessarily isoform (look at
the two black-filled elements).  

\end{slide} 

\begin{slide} 
\centerline{\tbf{\large Representation theorems}}

There are two types of representation theorems: \\
(a) The straight representation theorems.\\
(b) The \cpe results.\\

\end{slide} 
\begin{slide}
\centerline{\tbf{\large CPE defined}}

Let $K$ be a finite lattice.  

A finite lattice $L$ is a \emph{\cpe} of $K$, if $L$ is an
extension and every congruence $\gQ$  of $K$ has \emph{exactly one} extension
$\ol \gQ$ to~$L$---that is, $\ol \gQ|_K = \gQ$.  

Of course, then the congruence lattice of $K$ is isomorphic to the congruence lattice of
$L$. 

\bigskip

\centerline{\scalebox{.8}{\includegraphics{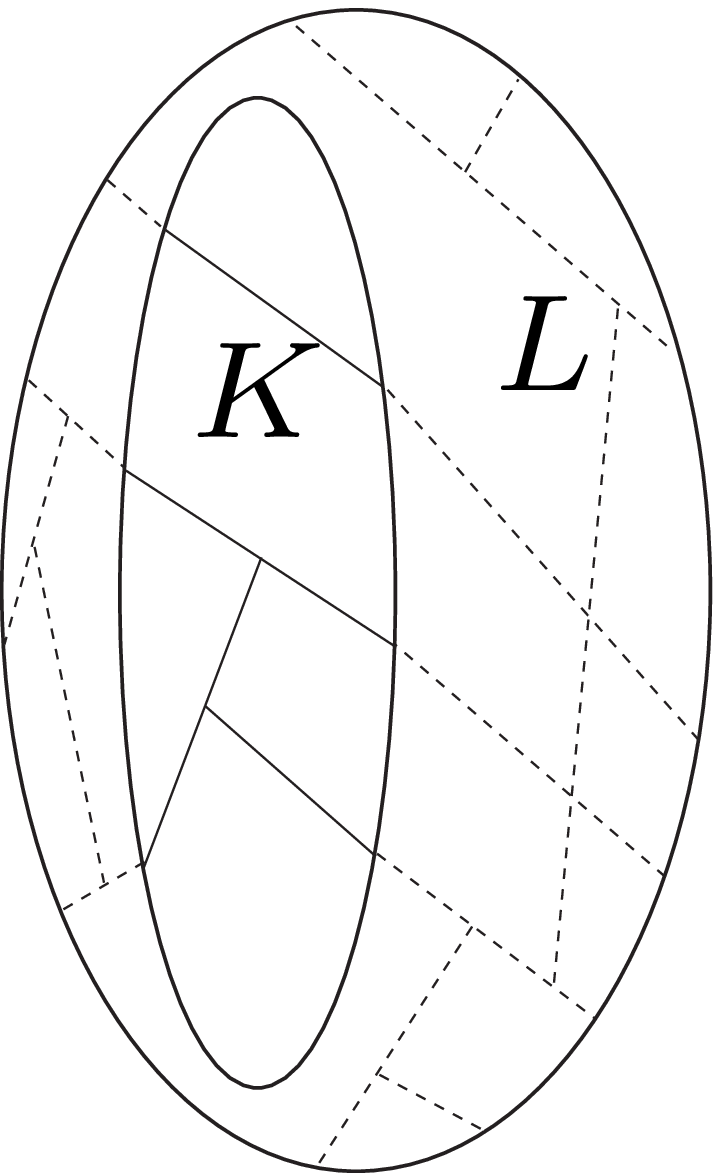}}}

\end{slide} 
\begin{slide}
\centerline{\tbf{\large CPE examples}}

\vspace{40pt}

\centerline{\includegraphics{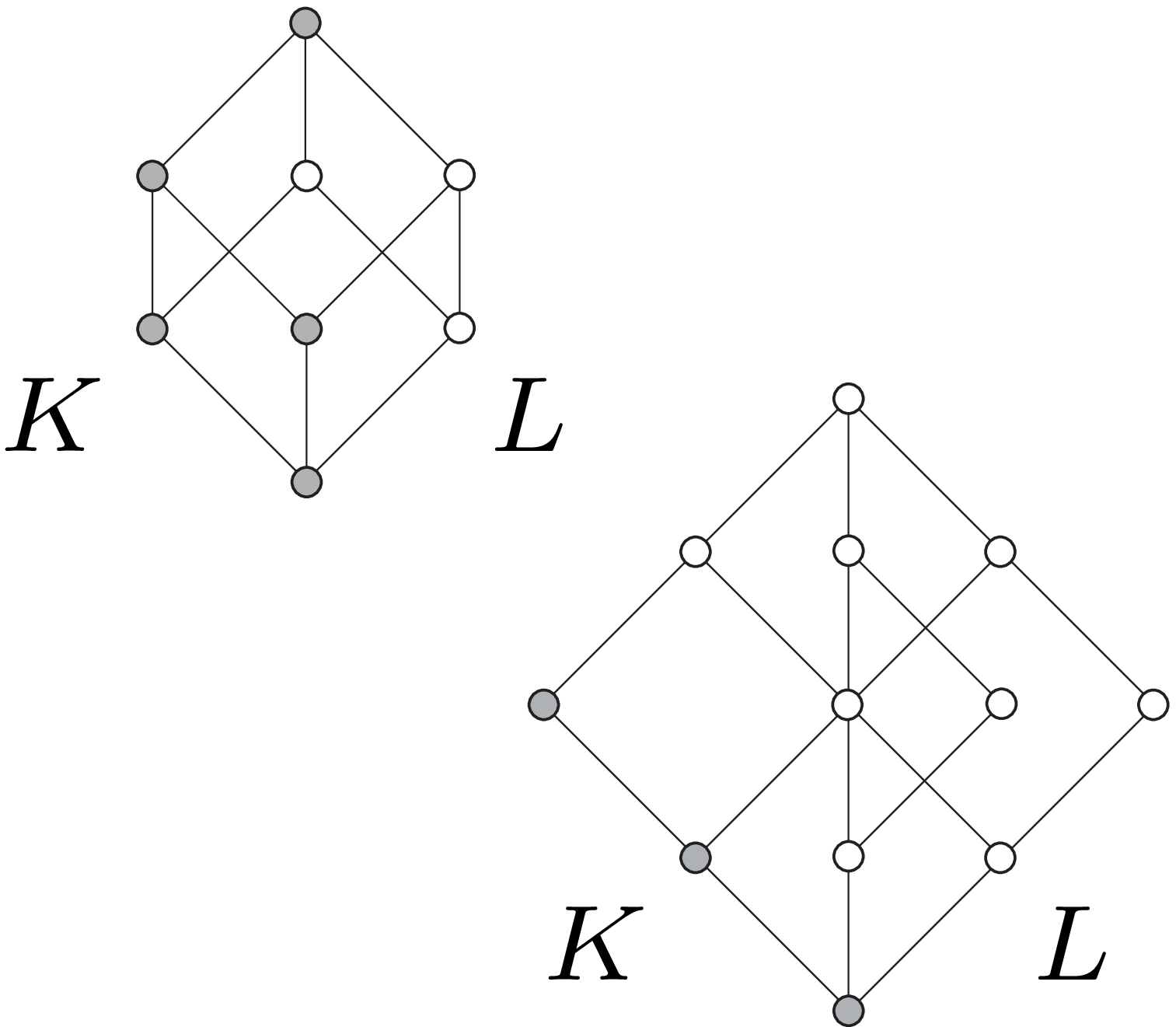}}

\end{slide} 
\begin{slide}
\centerline{\tbf{\large CPE fails}}

\vspace{40pt}

\centerline{\includegraphics{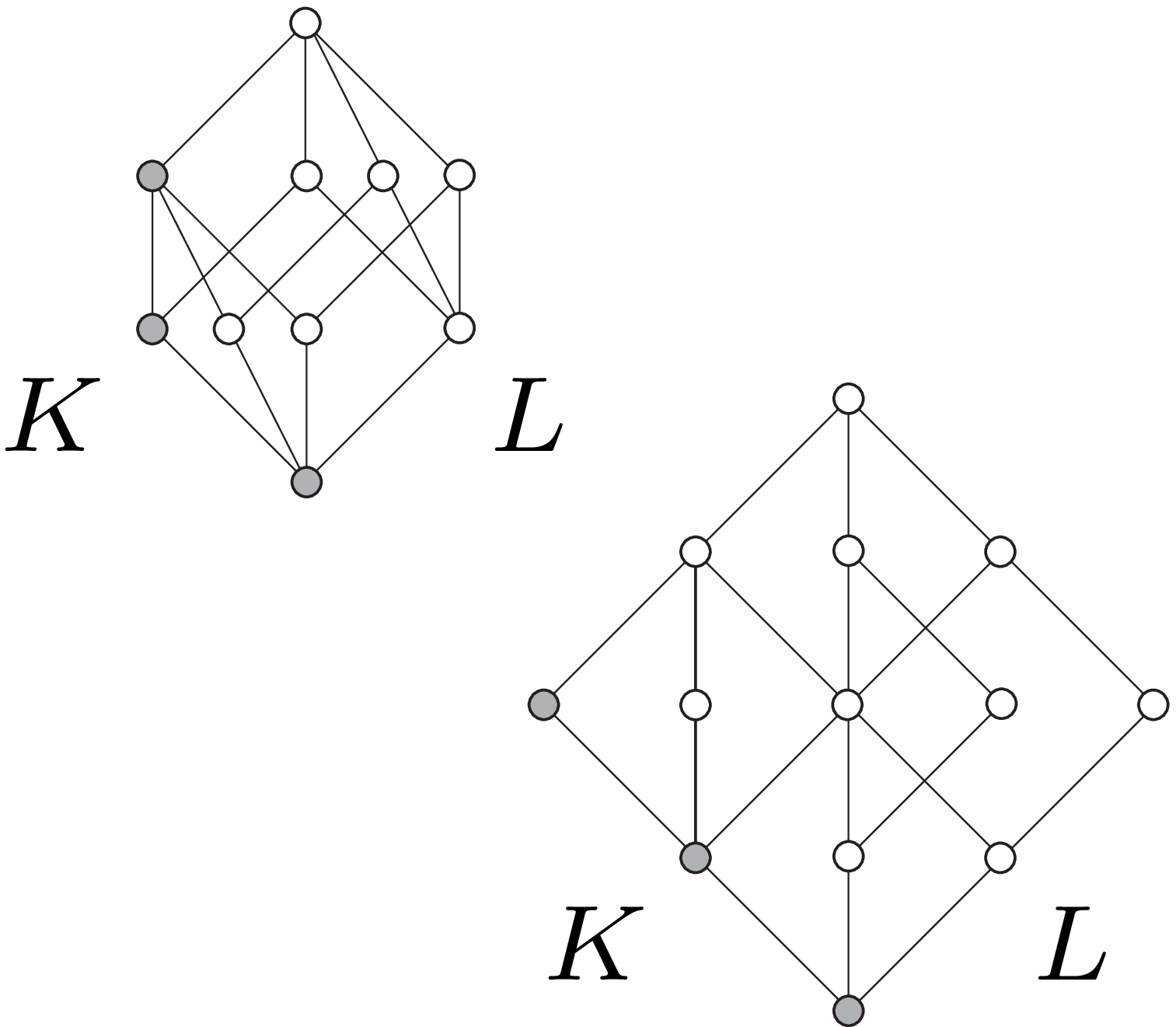}}

\end{slide} 

\begin{slide}

\centerline{\tbf{\large The basic technique for CPE results:}}
\centerline{\tbf{\large Cubic extensions}}

For a finite lattice $K$, let $\setm{K_i}{i \in I}$ be the subdirect factors of $K$.
For each $K_i$, we select $S(K_i)$, a finite \emph{simple} extension of $K_i$.

Let
\[
   \R(K) = \prodm{S(K_i)}{I \in I},
\]
and call it a \emph{cubic extension} of $K$. With the diagonal embedding, $K$ is a
sublattice of $\R(K)$.

The congruence lattice of $\R(K)$ is a finite Boolean
lattice; the number of atoms is the same as the number of join-irreducible
congruences of $K$. 

 Every congruence of $K$
extends to $\R(K)$.

\end{slide}

\begin{slide}
\centerline{\tbf{\large CPE result}}

G.~Gr\"atzer, E.\,T. Schmidt, and R. W. Quackenbush:

\bigskip

\begin{theorem}\label{T:isofoemcpe}
Every finite lattice $K$ has a \cpe to a finite isoform lattice $L$.
\end{theorem}

\end{slide} 
\begin{slide}
\centerline{\tbf{\large Separable lattices}}

A finite lattice $A$ is \emph{separable}, if it has an element $v$ that is a
\emph{separator}, that is, $0 \prec v \prec 1$.

\vspace{20pt}

\centerline{\includegraphics{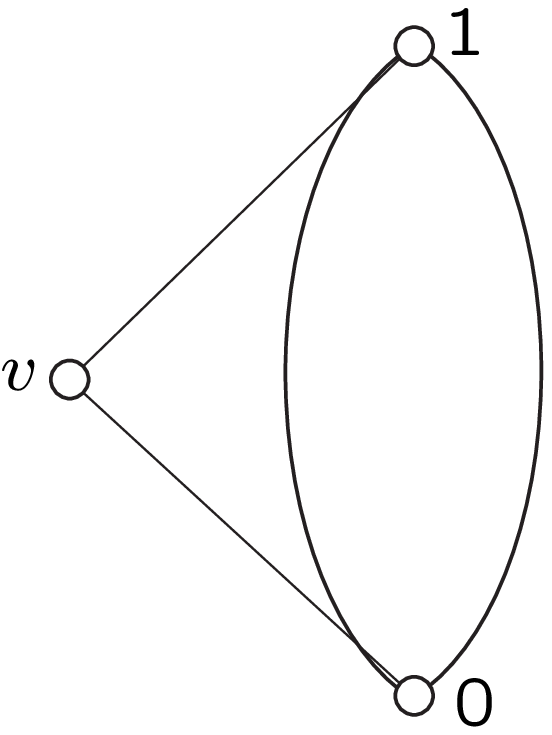}}

\end{slide}

\begin{slide}

\centerline{\tbf{\large Pruned lattices}}

For a finite poset $P$, and a family $S_p$, for $p \in P$, of
separable lattices, we construct a poset on the set $S = \prodm {S_p}{p \in P}$.

Let $\ba = \vv<a_p>_{p \in P},
\bb = \vv<b_p>_{p \in P} \in S$.

\begin{enumerate}
\item[(P)] $\ba \leq \bb$ in $L$ if{}f $\ba \leq_S \bb$ and if $p < p'$ in $P$, then
\[
   a_p =  v_p = b_p \text{\q implies that\q }a_{p'} = b_{p'}.
\]
\end{enumerate}
\end{slide} 

\begin{slide}

Examples:

If $P$ is unordered, then $\leq$ is the same as $\leq_S$. 

In the smallest
not unordered example, $P = \set{p, q}$ with $p < q$;
the figure below illustrates what we get: $S_p = S_q = C_3$,
the three-element chain; the two edges of $C_3^2$ that are not edges of $L$ are
dashed.

\centerline{\scalebox{1.5}{\includegraphics{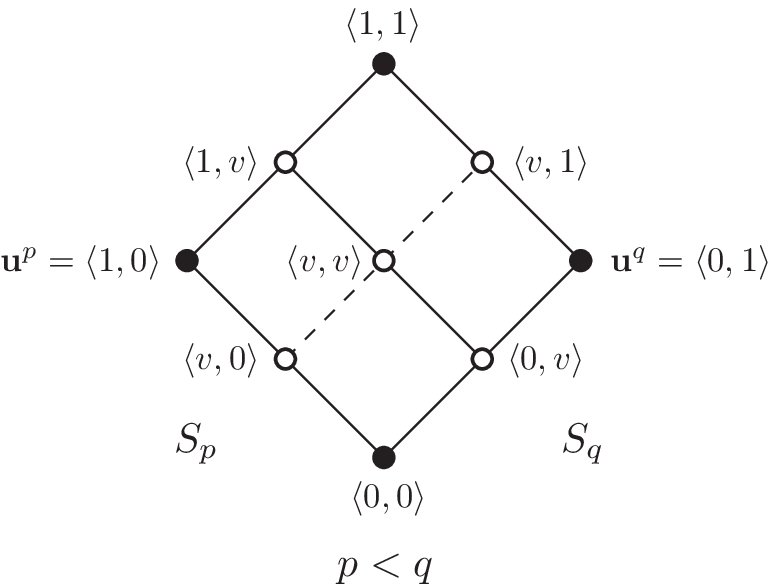}}}

\end{slide}
\begin{slide}

Next the poset ``hat'': $P = \set{p, q, r}$
with $p, r < q$, with the lattices $S_p = S_q = S_r= C_3 = \set{0, v, 1}$:

\centerline{\scalebox{1.4}{\includegraphics{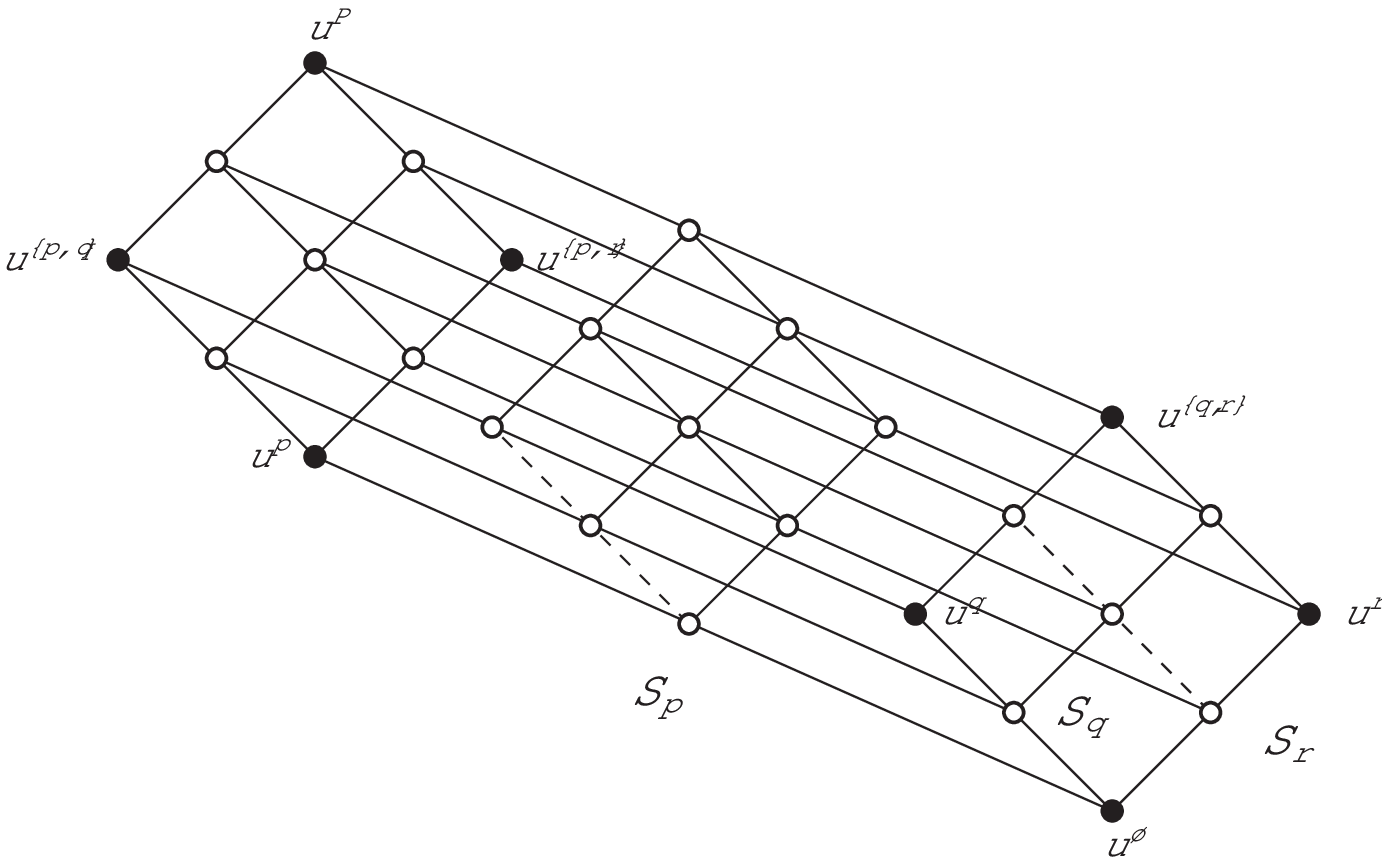}}}

Four edges of $C_3^3$ are missing in $L$; on the diagram these are marked with dashed
lines.

\end{slide} 

\begin{slide}
Applying this costruction, we have to assume that each $S_p$ is simple. So even the
smallest example is fairly large. Again, let $P = \set{p, q}$ with $p < q$ and $S_p = S_q =
M_3 = \set{0, a, b, v, 1}$. This is what we get:

\vspace{20pt}

\centerline{\scalebox{.7}{\includegraphics{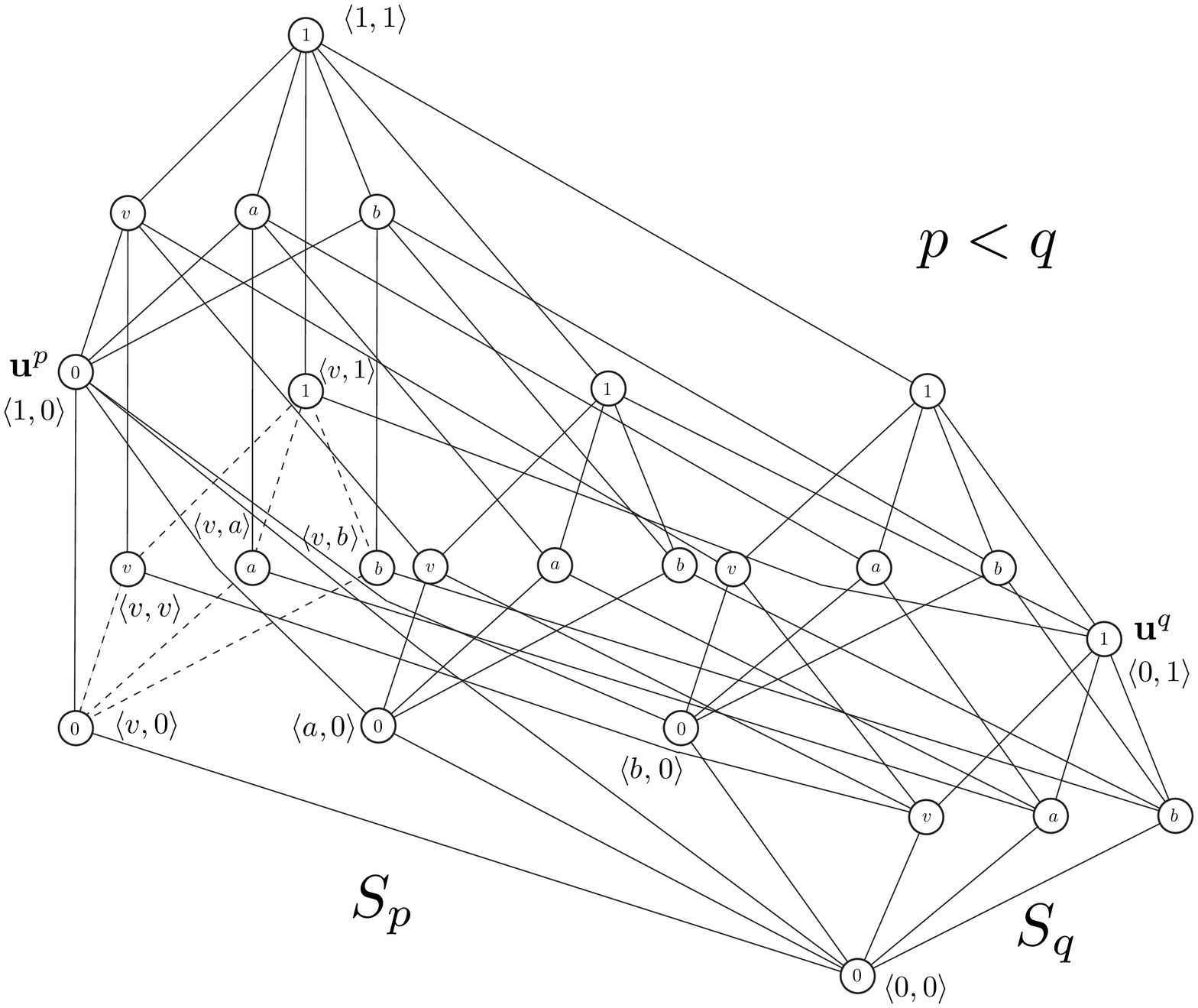}}}

\end{slide} 

\begin{slide}
In this diagram, $S_q$ is on the right, its elements are
labeled \text{$0$, $a$, $b$, $v$, $1$.} Its~unit element is $\bu^q = \vv<0, 1>$. 
The lattice $S_p$
is on the left, its elements are all labelled with $0$; its~unit element is 
$\bu^p = \vv<1, 0>$.
Five edges of $M_3^2$ are missing in $L$; on the diagram these are marked with 
dashed lines.

\end{slide} 

\begin{slide}
\centerline{\tbf{\large The fundamental theorem}}

Let $\ba = \vv<a_p>_{p \in P}$, $\bb = \vv<b_p>_{p
\in P} \in L$, and let $q \in P$; we shall call $q \in P$ an $\set{\ba,
\bb}$-\emph{fork}, if  $a_q = b_q = v$ and $a_{q'} \ne b_{q'}$, for some $q' > q$. 

\vspace{20pt}

\begin{theorem}\label{T:lattice}
$L$ is a lattice. Let $\ba = \vv<a_p>_{p \in P}$, $\bb = \vv<b_p>_{p \in P} \in L$.
Then
\[
(\ba \jj \bb)_p =
  \begin{cases}
            1, &\text{if $a_p \jj b_p = v$ and $\exists\, p' \geq p$,}\\
               &\text{\quad\tup{(1)}\, $p'$ is an $\set{\ba, \bb}$-fork, or}\\
               &\text{\quad\tup{(2)}\, $b_{p} \leq a_{p}$, $b_{p'} \nleq a_{p'}$, or}\\
               &\text{\quad\tup{(3)}\, $a_{p} \leq b_{p}$, $a_{p'} \nleq b_{p'}$;}\\
            a_p \jj b_p, &\text{otherwise};
   \end{cases}
\]
and
\[
(\ba \mm \bb)_p =
  \begin{cases}
            0, &\text{if $a_p \mm b_p = v$ and $\exists\, p' \geq p$,}\\
               &\text{\quad\tup{(1)}\, $p'$ is an $\set{\ba, \bb}$-fork, or}\\
               &\text{\quad\tup{(2)}\, $b_{p} \geq a_{p}$, $b_{p'} \ngeq a_{p'}$, or}\\
               &\text{\quad\tup{(3)}\, $a_{p} \geq b_{p}$, $a_{p'} \ngeq b_{p'}$;}\\
            a_p \mm b_p, &\text{otherwise}.
   \end{cases}
\]
\end{theorem}

\end{slide} 

\begin{slide}
Why not draw less trivial examples?

If we start with the lattice:

\vspace{20pt}

\centerline{\includegraphics{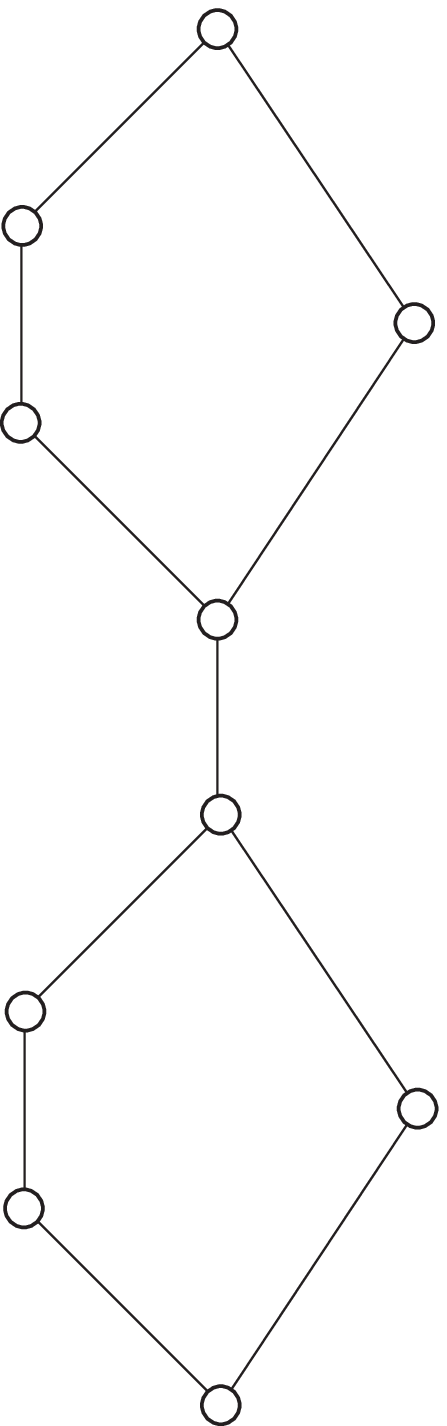}} 

then the uniform lattice we conctruct has 153,125 elements (and order dimension 7).
\end{slide} 
\begin{slide}
\centerline{\tbf{\large Problems}}

\vspace{20pt}

\begin{problem}
Are there analogues of these theorems for infinite lattices?
\end{problem}

\vspace{20pt}

\begin{problem}
Is every finite isoform lattice congruence permutable?
\end{problem}

\vspace{20pt}

\begin{problem}
Does every lattice have a \cpe to a congruence permutable lattice?
\end{problem}

Let $L$ be a finite lattice. A congruence $\gQ$ of $L$ is \emph{algebraically isoform} if, for
every $a \in L$, there is a unary algebraic function $p(x)$ that is an
isomorphism between $0/\gQ$ and $a/\gQ$. The lattice $L$ is algebraically isoform, if all
congruences are algebraically isoform.

\vspace{20pt}

\begin{problem}
Does every finite lattice has a \cpe to an algebraically isoform finite lattice?
\end{problem}

The congruence restriction problem.
\end{slide} 

\end{document}